\documentclass[a4paper,pdftex,12pt]{article}

\usepackage{tikz}
\usetikzlibrary{angles, quotes, babel}
\usetikzlibrary{decorations.pathreplacing}
\usetikzlibrary{spy}
\tikzset{SpyStyle/.style={spy using outlines={rectangle, magnification=15, width=2cm, height=2cm, connect spies, blue!70!black}}}

\usepackage{pgfplots}
\pgfplotsset{compat=1.13}
\usepgfplotslibrary{patchplots}
\definecolor{Farbe}{rgb}{0.8, 0.8, 0.8}

\usepackage{caption}
\usepackage{subcaption}

\usepackage{times}
\usepackage{amsmath}
\usepackage{amsfonts}
\usepackage{amssymb}

\usepackage{tabularx}
\usepackage[linkcolor=black,urlcolor=black,citecolor=black]{hyperref}
\usepackage[all]{hypcap}
\hypersetup{linktocpage,colorlinks=true,pdfborder={0 0 0}}

\usepackage{geometry}
\title{4-regular planar unit triangle graphs without additional triangles}

\author{
	Mike Winkler$^1$\quad Peter Dinkelacker$^2$\quad Stefan Vogel$^3$
}

\date{
	\small$^1$Fakult\"at f\"ur Mathematik, Ruhr-Universit\"at Bochum, Germany, mike.winkler@ruhr-uni-bochum.de\\[2mm]
	$^2$Togostr. 79, 13351 Berlin, Germany, peter@grity.de\\[2mm]
	$^3$Raun, Dorfstr. 7, 08648 Bad Brambach, Germany, backebackekuchen16@gmail.com
}

\begin{document}
  
  \maketitle
  
  \begin{abstract} 
    In this article we proof the existence of 4-regular planar unit-distance graphs consisting only of unit triangles without additional triangles. It is shown that the smallest number of unit triangles is $\leq6422$.
  \end{abstract}
  
  \section{\large{Introduction}}
  
  In 1991 it was shown by Harborth\cite{Harborth} that the smallest number of non-overlapping vertex-to-vertex unit triangles in the plane is $\leq42$ in general (see Figure 1), and $\leq3800$ if additional triangles are not allowed (see Figure 2). The smallest example of an additional triangle is an equilateral triangle consisting of three unit triangles. The graphs mentioned by Harborth are not planar graphs, because vertices lie exactly on other edges. So it remained an open question, whether 4-regular \textit{planar} unit-distance graphs consisting only of unit triangles without additional triangles exists. In this paper we prove the existence of such graphs. The currently smallest known example consists of 6422 unit triangles and has a rotational symmetry of order 192.
  
  \newpage
  
  Figure 1 shows the only known examples of 4-regular planar unit-distance graphs with 42 unit triangles. Both graphs are rigid and constructed from six copies of the same rigid subgraph (gray). The equilateral triangle $(A,B,C)$ shows one of the additional triangles.
  
  \begin{figure}[h]
  	\begin{subfigure}[b]{0.5\linewidth}
  		\centering
  		
  		\xdef\LstPN{0}
  		\newif\ifDupe
  		\pgfplotsset{avoid dupes/.code={\Dupefalse
  				\xdef\anker{\DefaultTextposition} 
  				\foreach \X in \LstPN
  				{\pgfmathtruncatemacro{\itest}{ifthenelse(\X==\punktnummer,1,0)}
  					\ifnum\itest=1
  					\global\Dupetrue
  					\breakforeach
  					\fi}
  				\ifDupe
  				\typeout{\punktnummer\space ist\space ein\space Duplikat!}%
  				\xdef\punktnummer{} 
  				\else
  				\xdef\LstPN{\LstPN,\punktnummer}
  				\typeout{\punktnummer\space ist\space neu\space mit\space urprgl.\space Anker=\anker}
  				\foreach \X in \LstExcept
  				{\ifnum\X=\punktnummer
  					\xdef\anker{\AusnahmeTextposition}
  					\fi}
  				\typeout{\punktnummer\space ist\space neu\space mit\space Anker=\anker}
  				\fi}}
  			
  		\def\DefaultTextposition{east}
  		\def\AusnahmeTextposition{west}
  		\def\AusnahmeListe{}
  		\xdef\BeliebigesVorhandenesKoordinatenpaar{{2.65, 0.35 }} 
  		\colorlet{Kantenfarbe}{gray}
  		\colorlet{Punktfarbe}{black}
  		\def\Beschriftung{} 
  		\pgfplotsset{
  			x=5.5mm, y=5.5mm,  
  		}
  	
  			\xdef\LstExcept{\AusnahmeListe}
  			\pgfdeclarelayer{bg}    
  			\pgfsetlayers{bg,main}  
  			
  			\pgfmathsetmacro{\xAlias}{\BeliebigesVorhandenesKoordinatenpaar[0]}
  			\pgfmathsetmacro{\yAlias}{\BeliebigesVorhandenesKoordinatenpaar[1]}
  			
  			\begin{tikzpicture}
  			\begin{axis}[hide axis, 
  			colormap={kantenfarbe}{color=(Kantenfarbe) color=(Kantenfarbe)},
  			thin, 
  			]
  			\addplot+[mark size=0.5pt, 
  			mark options={Punktfarbe}, 
  			table/row sep=newline, 
  			patch, 
  			patch type=polygon,
  			vertex count=2, 
  			%
  			patch table with point meta={
  				Startpkt Endpkt colordata  \\
  				1 9 \\
  				1 10 \\
  				2 1 \\
  				3 1 \\
  				3 2 \\
  				4 3 \\
  				4 2 \\
  				5 4 \\
  				5 2 \\
  				5 61 \\
  				5 63 \\
  				6 3 \\
  				6 4 \\
  				6 7 \\
  				7 7 \\
  				8 7 \\
  				9 7 \\
  				9 8 \\
  				10 8 \\
  				10 9 \\
  				11 8 \\
  				11 10 \\
  				11 14 \\
  				11 16 \\
  				12 7 \\
  				12 6 \\
  				12 17 \\
  				12 18 \\
  				13 20 \\
  				13 21 \\
  				14 13 \\
  				15 13 \\
  				15 14 \\
  				16 14 \\
  				16 15 \\
  				17 15 \\
  				17 16 \\
  				17 18 \\
  				18 18 \\
  				19 18 \\
  				20 18 \\
  				20 19 \\
  				21 19 \\
  				21 20 \\
  				22 19 \\
  				22 21 \\
  				22 25 \\
  				22 27 \\
  				23 41 \\
  				23 43 \\
  				23 45 \\
  				23 47 \\
  				24 31 \\
  				24 32 \\
  				25 24 \\
  				26 24 \\
  				26 25 \\
  				27 25 \\
  				27 26 \\
  				28 26 \\
  				28 27 \\
  				28 29 \\
  				29 29 \\
  				30 29 \\
  				31 29 \\
  				31 30 \\
  				32 30 \\
  				32 31 \\
  				33 30 \\
  				33 32 \\
  				33 36 \\
  				33 38 \\
  				34 28 \\
  				34 29 \\
  				34 39 \\
  				34 40 \\
  				35 42 \\
  				35 43 \\
  				36 35 \\
  				37 35 \\
  				37 36 \\
  				38 36 \\
  				38 37 \\
  				39 37 \\
  				39 38 \\
  				39 40 \\
  				40 40 \\
  				41 40 \\
  				42 40 \\
  				42 41 \\
  				43 41 \\
  				43 42 \\
  				44 51 \\
  				44 52 \\
  				45 44 \\
  				46 44 \\
  				46 45 \\
  				47 45 \\
  				47 46 \\
  				48 46 \\
  				48 47 \\
  				48 49 \\
  				49 49 \\
  				50 49 \\
  				51 49 \\
  				51 50 \\
  				52 50 \\
  				52 51 \\
  				53 50 \\
  				53 52 \\
  				53 56 \\
  				53 58 \\
  				54 48 \\
  				54 49 \\
  				54 59 \\
  				54 60 \\
  				55 62 \\
  				55 63 \\
  				56 55 \\
  				57 55 \\
  				57 56 \\
  				58 56 \\
  				58 57 \\
  				59 57 \\
  				59 58 \\
  				59 60 \\
  				60 60 \\
  				61 60 \\
  				62 60 \\
  				62 61 \\
  				63 61 \\
  				63 62 \\
  			},
  			%
  			visualization depends on={value \thisrowno{0} \as \punktnummer},
  			every node near coord/.append style={
  				/pgfplots/avoid dupes,
  			},
  			nodes near coords={\Beschriftung},
  			nodes near coords style={
  				anchor=\anker,
  				text=black, font=\scriptsize, 
  				name=p-\punktnummer, 
  				path picture={
  					\coordinate[] (P\punktnummer) at (p-\punktnummer.\anker);}
  			},
  			]
  			table[header=true, x index=1, y index=2, row sep=\\] {
  				Nr x y        \\
  				0 0 0         \\
  				1 2.65 0.35   \\
  				2 3.65 0.28   \\
  				3 3.21 1.18   \\
  				4 4.21 1.11   \\
  				5 4.64 0.21   \\
  				6 3.77 2.01   \\
  				7 2.83 2.34   \\
  				8 1.92 1.92   \\
  				9 2.74 1.35   \\
  				10 1.83 0.93   \\
  				11 1.01 1.50   \\
  				12 3.59 2.99   \\
  				13 0.83 3.50   \\
  				14 0.92 2.50   \\
  				15 1.74 3.08   \\
  				16 1.83 2.08   \\
  				17 2.65 2.66   \\
  				18 2.83 3.64   \\
  				19 2.27 4.47   \\
  				20 1.83 3.57   \\
  				21 1.27 4.40   \\
  				22 1.71 5.30   \\
  				23 7.58 5.30   \\
  				24 2.83 6.95   \\
  				25 2.27 6.12   \\
  				26 3.26 6.05   \\
  				27 2.70 5.22   \\
  				28 3.70 5.15   \\
  				29 4.46 5.80   \\
  				30 4.55 6.80   \\
  				31 3.65 6.38   \\
  				32 3.74 7.37   \\
  				33 4.64 7.79   \\
  				34 4.64 4.82   \\
  				35 6.46 6.95   \\
  				36 5.55 7.37   \\
  				37 5.64 6.38   \\
  				38 4.73 6.80   \\
  				39 4.82 5.80   \\
  				40 5.58 5.15   \\
  				41 6.58 5.22   \\
  				42 6.02 6.05   \\
  				43 7.02 6.12   \\
  				44 8.45 3.50   \\
  				45 8.02 4.40   \\
  				46 7.46 3.57   \\
  				47 7.02 4.47   \\
  				48 6.46 3.64   \\
  				49 6.64 2.66   \\
  				50 7.46 2.08   \\
  				51 7.55 3.08   \\
  				52 8.36 2.50   \\
  				53 8.27 1.50   \\
  				54 5.70 2.99   \\
  				55 6.64 0.35   \\
  				56 7.46 0.93   \\
  				57 6.55 1.35   \\
  				58 7.37 1.92   \\
  				59 6.46 2.34   \\
  				60 5.52 2.01   \\
  				61 5.08 1.11   \\
  				62 6.08 1.18   \\
  				63 5.64 0.28   \\
  			};
  			
  			\addplot[no marks, 
  			nodes near coords={},
  			visualization depends on={value \thisrowno{0} \as \PunktI},
  			visualization depends on={value \thisrowno{1} \as \Scheitel},
  			visualization depends on={value \thisrowno{2} \as \PunktII},
  			visualization depends on={value \thisrowno{3} \as \Winkelradius},
  			visualization depends on={value \thisrowno{4} \as \Winkelfarbe},
  			nodes near coords style={anchor=center,
  				path picture={
  					\begin{pgfonlayer}{bg}    
  					\draw pic [angle radius=\Winkelradius cm,
  					fill=\Winkelfarbe!40, draw=\Winkelfarbe, 
  					] {angle = P\PunktI--P\Scheitel--P\PunktII};
  					\end{pgfonlayer}
  			}},%
  			]
  			table[header=true, x expr =2.65, y expr=0.35]{
  				Punkt1 Scheitel Punkt2 Winkelradius[cm] Winkelfarbe 
  			};
  			
  			\end{axis}
  			
  			\begin{pgfonlayer}{bg}
  			\fill[Farbe] (P13) -- (P20) -- (P21) -- cycle;
  			\fill[Farbe] (P20) -- (P18) -- (P19) -- cycle;
  			\fill[Farbe] (P21) -- (P19) -- (P22) -- cycle;
  			\fill[Farbe] (P13) -- (P14) -- (P15) -- cycle;
  			\fill[Farbe] (P15) -- (P16) -- (P17) -- cycle;
  			\fill[Farbe] (P11) -- (P14) -- (P16) -- cycle;
  			\fill[Farbe] (P12) -- (P17) -- (P18) -- cycle;
  			\node[left]  at (P22) {$A$};
  			\node[left]  at (P13) {$B$};
  			\node[below] at (P19) {$\quad\quad\quad C$};
  			\end{pgfonlayer}
  			
  			\end{tikzpicture}
  	\end{subfigure}
  	\begin{subfigure}[b]{0.4\linewidth}
  		\centering
  		
  		\xdef\LstPN{0}
  		\newif\ifDupe
  		\pgfplotsset{avoid dupes/.code={\Dupefalse
  				\xdef\anker{\DefaultTextposition} 
  				\foreach \X in \LstPN
  				{\pgfmathtruncatemacro{\itest}{ifthenelse(\X==\punktnummer,1,0)}
  					\ifnum\itest=1
  					\global\Dupetrue
  					\breakforeach
  					\fi}
  				\ifDupe
  				\typeout{\punktnummer\space ist\space ein\space Duplikat!}%
  				\xdef\punktnummer{} 
  				\else
  				\xdef\LstPN{\LstPN,\punktnummer}
  				\typeout{\punktnummer\space ist\space neu\space mit\space urprgl.\space Anker=\anker}
  				\foreach \X in \LstExcept
  				{\ifnum\X=\punktnummer
  					\xdef\anker{\AusnahmeTextposition}
  					\fi}
  				\typeout{\punktnummer\space ist\space neu\space mit\space Anker=\anker}
  				\fi}}
  			
  		\def\DefaultTextposition{east}
  		\def\AusnahmeTextposition{west}
  		\def\AusnahmeListe{}
  		\xdef\BeliebigesVorhandenesKoordinatenpaar{{2.40, 0.20 }} 
  		\colorlet{Kantenfarbe}{gray}
  		\colorlet{Punktfarbe}{black}
  		\def\Beschriftung{} 
  		\pgfplotsset{
  			x=5.5mm, y=5.5mm,  
  		}
  		
  			\xdef\LstExcept{\AusnahmeListe}
  			\pgfdeclarelayer{bg}    
  			\pgfsetlayers{bg,main}  
  			
  			\pgfmathsetmacro{\xAlias}{\BeliebigesVorhandenesKoordinatenpaar[0]}
  			\pgfmathsetmacro{\yAlias}{\BeliebigesVorhandenesKoordinatenpaar[1]}
  			
  			\begin{tikzpicture} 
  			\begin{axis}[hide axis, 
  			colormap={kantenfarbe}{color=(Kantenfarbe) color=(Kantenfarbe)},
  			thin, 
  			]
  			\addplot+[mark size=0.5pt, 
  			mark options={Punktfarbe}, 
  			table/row sep=newline, 
  			patch, 
  			patch type=polygon,
  			vertex count=2, 
  			%
  			patch table with point meta={
  				Startpkt Endpkt colordata  \\
  				1 9 \\
  				1 10 \\
  				2 1 \\
  				3 1 \\
  				3 2 \\
  				4 3 \\
  				4 2 \\
  				5 4 \\
  				5 2 \\
  				5 43 \\
  				6 3 \\
  				6 4 \\
  				6 7 \\
  				7 7 \\
  				8 7 \\
  				9 7 \\
  				9 8 \\
  				10 8 \\
  				10 9 \\
  				11 8 \\
  				11 10 \\
  				11 14 \\
  				11 16 \\
  				12 7 \\
  				12 6 \\
  				12 17 \\
  				12 18 \\
  				13 20 \\
  				13 21 \\
  				14 13 \\
  				15 13 \\
  				15 14 \\
  				16 14 \\
  				16 15 \\
  				17 15 \\
  				17 16 \\
  				17 18 \\
  				18 18 \\
  				19 18 \\
  				20 18 \\
  				20 19 \\
  				21 19 \\
  				21 20 \\
  				22 19 \\
  				22 21 \\
  				22 46 \\
  				22 48 \\
  				23 5 \\
  				23 41 \\
  				24 32 \\
  				24 33 \\
  				25 24 \\
  				26 24 \\
  				26 25 \\
  				27 25 \\
  				27 26 \\
  				28 25 \\
  				28 27 \\
  				28 55 \\
  				28 57 \\
  				29 26 \\
  				29 27 \\
  				29 30 \\
  				30 30 \\
  				31 30 \\
  				32 30 \\
  				32 31 \\
  				33 31 \\
  				33 32 \\
  				34 31 \\
  				34 33 \\
  				34 37 \\
  				34 39 \\
  				35 29 \\
  				35 30 \\
  				35 40 \\
  				35 41 \\
  				36 42 \\
  				36 43 \\
  				37 36 \\
  				38 36 \\
  				38 37 \\
  				39 37 \\
  				39 38 \\
  				40 38 \\
  				40 39 \\
  				40 41 \\
  				41 41 \\
  				42 41 \\
  				42 23 \\
  				43 23 \\
  				43 42 \\
  				44 51 \\
  				44 53 \\
  				44 58 \\
  				44 59 \\
  				45 52 \\
  				45 53 \\
  				46 45 \\
  				47 45 \\
  				47 46 \\
  				48 46 \\
  				48 47 \\
  				49 47 \\
  				49 48 \\
  				49 50 \\
  				50 50 \\
  				51 50 \\
  				52 50 \\
  				52 51 \\
  				53 51 \\
  				53 52 \\
  				54 61 \\
  				54 62 \\
  				55 54 \\
  				56 54 \\
  				56 55 \\
  				57 55 \\
  				57 56 \\
  				58 56 \\
  				58 57 \\
  				58 59 \\
  				59 59 \\
  				60 59 \\
  				61 59 \\
  				61 60 \\
  				62 60 \\
  				62 61 \\
  				63 60 \\
  				63 62 \\
  				63 49 \\
  				63 50 \\
  			},
  			%
  			visualization depends on={value \thisrowno{0} \as \punktnummer},
  			every node near coord/.append style={
  				/pgfplots/avoid dupes,
  			},
  			nodes near coords={\Beschriftung},
  			nodes near coords style={
  				anchor=\anker,
  				text=black, font=\scriptsize, 
  				name=p-\punktnummer, 
  				path picture={
  					\coordinate[] (P\punktnummer) at (p-\punktnummer.\anker);}
  			},
  			]
  			table[header=true, x index=1, y index=2, row sep=\\] {
  				Nr x y        \\
  				0 0 0         \\
  				1 2.40 0.20   \\
  				2 3.40 0.16   \\
  				3 2.93 1.05   \\
  				4 3.93 1.01   \\
  				5 4.40 0.12   \\
  				6 3.47 1.89   \\
  				7 2.51 2.20   \\
  				8 1.62 1.75   \\
  				9 2.46 1.20   \\
  				10 1.56 0.75   \\
  				11 0.73 1.30   \\
  				12 3.25 2.87   \\
  				13 0.48 3.28   \\
  				14 0.60 2.29   \\
  				15 1.40 2.89   \\
  				16 1.53 1.90   \\
  				17 2.32 2.50   \\
  				18 2.47 3.49   \\
  				19 1.88 4.30   \\
  				20 1.48 3.39   \\
  				21 0.89 4.19   \\
  				22 1.29 5.11   \\
  				23 4.86 1.01   \\
  				24 8.32 3.28   \\
  				25 7.91 4.19   \\
  				26 7.32 3.39   \\
  				27 6.92 4.30   \\
  				28 7.50 5.11   \\
  				29 6.33 3.49   \\
  				30 6.47 2.50   \\
  				31 7.27 1.90   \\
  				32 7.39 2.89   \\
  				33 8.19 2.29   \\
  				34 8.07 1.30   \\
  				35 5.54 2.87   \\
  				36 6.40 0.20   \\
  				37 7.23 0.75   \\
  				38 6.34 1.20   \\
  				39 7.18 1.75   \\
  				40 6.28 2.20   \\
  				41 5.33 1.89   \\
  				42 5.86 1.05   \\
  				43 5.40 0.16   \\
  				44 4.88 6.52   \\
  				45 2.89 6.31   \\
  				46 2.09 5.71   \\
  				47 3.01 5.32   \\
  				48 2.21 4.72   \\
  				49 3.13 4.33   \\
  				50 4.06 4.69   \\
  				51 4.47 5.60   \\
  				52 3.48 5.50   \\
  				53 3.89 6.41   \\
  				54 5.90 3.91   \\
  				55 6.70 4.51   \\
  				56 5.78 4.90   \\
  				57 6.58 5.50   \\
  				58 5.66 5.89   \\
  				59 4.73 5.53   \\
  				60 4.32 4.61   \\
  				61 5.32 4.72   \\
  				62 4.91 3.80   \\
  				63 3.92 3.70   \\
  			};
  			
  			\addplot[no marks, 
  			nodes near coords={},
  			visualization depends on={value \thisrowno{0} \as \PunktI},
  			visualization depends on={value \thisrowno{1} \as \Scheitel},
  			visualization depends on={value \thisrowno{2} \as \PunktII},
  			visualization depends on={value \thisrowno{3} \as \Winkelradius},
  			visualization depends on={value \thisrowno{4} \as \Winkelfarbe},
  			nodes near coords style={anchor=center,
  				path picture={
  					\begin{pgfonlayer}{bg}    
  					\draw pic [angle radius=\Winkelradius cm,
  					fill=\Winkelfarbe!40, draw=\Winkelfarbe, 
  					] {angle = P\PunktI--P\Scheitel--P\PunktII};
  					\end{pgfonlayer}
  			}},%
  			]
  			table[header=true, x expr =2.40, y expr=0.20]{
  				Punkt1 Scheitel Punkt2 Winkelradius[cm] Winkelfarbe 
  			};
  			
  			\end{axis}
  			
  			
  			\end{tikzpicture}
  	\end{subfigure}
  	\caption{}
  \end{figure}
  
  For graphs without additional triangles Harborth presented the subgraph $G_1$ which consists of 38 unit triangles (see Figure 2).
  
  \begin{center}
  	\begin{minipage}{\linewidth}
  		\centering
  		
  		\xdef\LstPN{0}
  		\newif\ifDupe
  		\pgfplotsset{avoid dupes/.code={\Dupefalse
  				\xdef\anker{\DefaultTextposition} 
  				\foreach \X in \LstPN
  				{\pgfmathtruncatemacro{\itest}{ifthenelse(\X==\punktnummer,1,0)}
  					\ifnum\itest=1
  					\global\Dupetrue
  					\breakforeach
  					\fi}
  				\ifDupe
  				\typeout{\punktnummer\space ist\space ein\space Duplikat!}%
  				\xdef\punktnummer{} 
  				\else
  				\xdef\LstPN{\LstPN,\punktnummer}
  				\typeout{\punktnummer\space ist\space neu\space mit\space urprgl.\space Anker=\anker}
  				\foreach \X in \LstExcept
  				{\ifnum\X=\punktnummer
  					\xdef\anker{\AusnahmeTextposition}
  					\fi}
  				\typeout{\punktnummer\space ist\space neu\space mit\space Anker=\anker}
  				\fi}}
  		
  		\def\DefaultTextposition{east}
  		\def\AusnahmeTextposition{west}
  		\def\AusnahmeListe{}
  		\xdef\BeliebigesVorhandenesKoordinatenpaar{{0.86645437967485794406, 2.46197022175259672139}} 
  		\colorlet{Kantenfarbe}{gray}
  		\colorlet{Punktfarbe}{black}
  		\def\Beschriftung{} 
  		\pgfplotsset{
  			x=12.0mm, y=12.0mm,  
  		}
  		
  			\xdef\LstExcept{\AusnahmeListe}
  			\pgfdeclarelayer{bg}    
  			\pgfsetlayers{bg,main}  
  			
  			\pgfmathsetmacro{\xAlias}{\BeliebigesVorhandenesKoordinatenpaar[0]}
  			\pgfmathsetmacro{\yAlias}{\BeliebigesVorhandenesKoordinatenpaar[1]}
  			
  			\begin{tikzpicture}[SpyStyle]
  			
  			
  			\begin{axis}[hide axis, colormap={kantenfarbe}{color=(Kantenfarbe) color=(Kantenfarbe)}, line width=0.1, 
  			]
  			\addplot+[mark size=0.01pt, 
  			mark options={Punktfarbe}, 
  			table/row sep=newline, 
  			patch, 
  			patch type=polygon,
  			vertex count=2, 
  			%
  			patch table with point meta={
  				Startpkt Endpkt colordata  \\
  				1 1 \\
  				2 1 \\
  				3 1 \\
  				3 2 \\
  				3 39 \\
  				3 40 \\
  				4 1 \\
  				5 1 \\
  				5 4 \\
  				6 5 \\
  				7 5 \\
  				7 6 \\
  				8 7 \\
  				9 7 \\
  				9 8 \\
  				10 9 \\
  				11 9 \\
  				11 10 \\
  				12 11 \\
  				13 11 \\
  				13 12 \\
  				14 13 \\
  				15 13 \\
  				15 14 \\
  				16 15 \\
  				17 15 \\
  				17 16 \\
  				18 17 \\
  				19 17 \\
  				19 18 \\
  				20 19 \\
  				20 56 \\
  				21 19 \\
  				21 20 \\
  				22 2 \\
  				23 2 \\
  				23 22 \\
  				24 22 \\
  				25 24 \\
  				25 22 \\
  				26 25 \\
  				27 26 \\
  				27 25 \\
  				28 27 \\
  				29 28 \\
  				29 27 \\
  				30 29 \\
  				31 30 \\
  				31 29 \\
  				32 31 \\
  				33 32 \\
  				33 31 \\
  				34 33 \\
  				35 34 \\
  				35 33 \\
  				36 35 \\
  				37 36 \\
  				37 35 \\
  				37 21 \\
  				38 37 \\
  				38 21 \\
  				39 39 \\
  				40 39 \\
  				41 39 \\
  				42 39 \\
  				42 41 \\
  				43 42 \\
  				44 42 \\
  				44 43 \\
  				45 44 \\
  				46 44 \\
  				46 45 \\
  				47 46 \\
  				48 46 \\
  				48 47 \\
  				49 48 \\
  				50 48 \\
  				50 49 \\
  				51 50 \\
  				52 50 \\
  				52 51 \\
  				53 52 \\
  				54 52 \\
  				54 53 \\
  				55 54 \\
  				56 54 \\
  				56 55 \\
  				57 56 \\
  				57 20 \\
  				58 40 \\
  				59 40 \\
  				59 58 \\
  				60 58 \\
  				61 58 \\
  				61 60 \\
  				62 61 \\
  				63 61 \\
  				63 62 \\
  				64 63 \\
  				65 63 \\
  				65 64 \\
  				66 65 \\
  				67 65 \\
  				67 66 \\
  				68 67 \\
  				69 67 \\
  				69 68 \\
  				70 69 \\
  				71 69 \\
  				71 70 \\
  				72 71 \\
  				73 57 \\
  				73 71 \\
  				73 72 \\
  				74 57 \\
  				74 73 \\
  			},
  			%
  			visualization depends on={value \thisrowno{0} \as \punktnummer},
  			every node near coord/.append style={
  				/pgfplots/avoid dupes,
  			},
  			nodes near coords={\Beschriftung},
  			nodes near coords style={
  				anchor=\anker,
  				text=black, font=\scriptsize, 
  				name=p-\punktnummer, 
  				path picture={
  					\coordinate[] (P\punktnummer) at (p-\punktnummer.\anker);}
  			},
  			]
  			table[header=true, x index=1, y index=2, row sep=\\] {
  				Nr x y        \\
  				0 0 0         \\
  				1 0.86645437967485794406 2.46197022175259672139  \\
  				2 0.00000000000000000000 2.96122647653790949107  \\
  				3 0.00085859019507165913 1.96122684512653933098  \\
  				4 1.73205016915464637961 1.96122684512654266165  \\
  				5 1.73290875934971344563 2.96122647653791304379  \\
  				6 1.73376734954478517459 1.96122684512654266165  \\
  				7 2.59936313902457083458 2.46197022175260071819  \\
  				8 3.46495892850435982524 1.96122684512654754663  \\
  				9 3.46581751869942689126 2.96122647653791792877  \\
  				10 3.46667610889449839817 1.96122684512654754663  \\
  				11 4.33227189837428472430 2.46197022175260471499  \\
  				12 5.19786768785407371496 1.96122684512655087730  \\
  				13 5.19872627804914078098 2.96122647653792192557  \\
  				14 5.19958486824421228789 1.96122684512655154343  \\
  				15 6.06518065772399861402 2.46197022175260871180  \\
  				16 6.93077644720378671650 1.96122684512655509614  \\
  				17 6.93163503739885378252 2.96122647653792547828  \\
  				18 6.93249362759392617761 1.96122684512655509614  \\
  				19 7.79808941707371250374 2.46197022175261359678  \\
  				20 8.66368520655350060622 1.96122684512656042521  \\
  				21 8.66454379674856767224 2.96122647653793080735  \\
  				22 0.97032597650811691636 3.20302701275230550237  \\
  				23 0.27575758124369248447 3.92245369025307910604  \\
  				24 1.70869697454706703788 3.87742175565512869895  \\
  				25 1.92355445506011202284 2.90077634248431603936  \\
  				26 2.19931203630380389669 3.86200355619948565433  \\
  				27 2.89388043156822893920 3.14257687869871249475  \\
  				28 3.63225142960717839458 3.81697162160153657950  \\
  				29 3.84710891012022404567 2.84032620843072391992  \\
  				30 4.12286649136391503134 3.80155342214589442307  \\
  				31 4.81743488662834096203 3.08212674464512126349  \\
  				32 5.55580588466729086150 3.75652148754794446006  \\
  				33 5.77066336518033562442 2.77987607437713180047  \\
  				34 6.04642094642402749827 3.74110328809230141545  \\
  				35 6.74098934168845165260 3.02167661059152825587  \\
  				36 7.47936033972740155207 3.69607135349435145244  \\
  				37 7.69421782024044631498 2.71942594032353923694  \\
  				38 7.96997540148413818883 3.68065315403870885191  \\
  				39 0.86645437967486060860 1.46048346850048615941  \\
  				40 0.00000000000000479872 0.96122721371516894884  \\
  				41 1.73205016915464660165 1.96122684512654399391  \\
  				42 1.73290875934971855266 0.96122721371517438893  \\
  				43 1.73376734954478561868 1.96122684512654399391  \\
  				44 2.59936313902457349911 1.46048346850049015622  \\
  				45 3.46495892850436026933 1.96122684512654887889  \\
  				46 3.46581751869943222033 0.96122721371517849676  \\
  				47 3.46667610889449928635 1.96122684512654887889  \\
  				48 4.33227189837428827701 1.46048346850049437506  \\
  				49 5.19786768785407460314 1.96122684512655354183  \\
  				50 5.19872627804914699823 0.96122721371518060618  \\
  				51 5.19958486824421317607 1.96122684512655354183  \\
  				52 6.06518065772399950220 1.46048346850049903800  \\
  				53 6.93077644720378760468 1.96122684512655842681  \\
  				54 6.93163503739885733523 0.96122721371518882183  \\
  				55 6.93249362759392706579 1.96122684512655842681  \\
  				56 7.79808941707371428009 1.46048346850050525525  \\
  				57 8.66454379674857122495 0.96122721371519015410  \\
  				58 0.97032597650812280055 0.71942667750077782252  \\
  				59 0.27575758124370192137 0.00000000000000000000  \\
  				60 1.70869697454707702988 0.04503193459795817866  \\
  				61 1.92355445506011646373 1.02167734776877128233  \\
  				62 2.19931203630381322256 0.06045013405360445680  \\
  				63 2.89388043156823560054 0.77987681155438026703  \\
  				64 3.63225142960718816454 0.10548206865155988765  \\
  				65 3.84710891012022893065 1.08212748182237228356  \\
  				66 4.12286649136392480131 0.12090026810720343187  \\
  				67 4.81743488662834717928 0.84032694560798160133  \\
  				68 5.55580588466729885511 0.16593220270516056969  \\
  				69 5.77066336518033828895 1.14257761587597417297  \\
  				70 6.04642094642403549187 0.18135040216080686171  \\
  				71 6.74098934168845698167 0.90077707966158127029  \\
  				72 7.47936033972741132203 0.22638233675876504036  \\
  				73 7.69421782024045075588 1.20302774992957983713  \\
  				74 7.96997540148414618244 0.24180053621440686373  \\
  			};
  			
  			\addplot[no marks, 
  			nodes near coords={},
  			visualization depends on={value \thisrowno{0} \as \PunktI},
  			visualization depends on={value \thisrowno{1} \as \Scheitel},
  			visualization depends on={value \thisrowno{2} \as \PunktII},
  			visualization depends on={value \thisrowno{3} \as \Winkelradius},
  			visualization depends on={value \thisrowno{4} \as \Winkelfarbe},
  			visualization depends on={value \thisrowno{5} \as \Winkelname},
  			visualization depends on={value \thisrowno{6} \as \WinkelExzentrizitaet},
  			nodes near coords style={anchor=center,
  				path picture={
  					\begin{pgfonlayer}{bg}    
  					\draw pic [angle radius=\Winkelradius cm,%
  					fill=\Winkelfarbe!40, draw=\Winkelfarbe,
  					"$\Winkelname$", angle eccentricity =\WinkelExzentrizitaet,
  					text=\Winkelfarbe%
  					] {angle = P\PunktI--P\Scheitel--P\PunktII};
  					\end{pgfonlayer}
  			}},%
  			]
  			table[header=true, x expr=0.86645437967485794406, y expr=2.46197022175259672139]{
  				Punkt1 Scheitel Punkt2 Winkelradius[cm] Winkelfarbe Winkelname WinkelExz
  				3 1 4 0.5 blue!70!black {\gamma} 1.4 \\
  				6 7 8 0.5 blue!70!black \gamma{} 1.4 \\
  				10 11 12 0.5 blue!70!black {\gamma} 1.4 \\
  				14 15 16 0.5 blue!70!black {\gamma} 1.4 \\
  				18 19 20 0.5 blue!70!black {\gamma} 1.4 \\
  				1 2 22 0.4 green!50!black {\beta} 1.4 \\
  				24 22 23 0.4 orange!80!black {\alpha} 1.4 \\
  				28 27 26 0.4 orange!80!black {\alpha} 1.4 \\
  				32 31 30 0.4 orange!80!black {\alpha} 1.4 \\
  				36 35 34 0.4 orange!80!black {\alpha} 1.4 \\
  			};
  			
  			\end{axis}
  			
  			\draw[blue!70!white, line width=0.1, shorten >= -35, shorten <=-30] (P23) -- (P38);
  			\draw[blue!70!white, line width=0.1, shorten >= -35, shorten <=-30] (P59) -- (P74);
  			\node[below, blue!70!white, font=\Large] at (P23) {$g_{1}\quad\quad\quad$};
  			\node[above, blue!70!white, font=\Large] at (P59) {$g_{2}\quad\quad\quad$};
  			
  			\node[above] at (P23) {$A$};
  			\node[left]  at (P3) {$B$};
  			\node[below] at (P59) {$C$};
  			\node[below] at (P74) {$D$};
  			\node[right] at (P20) {$E$};
  			\node[above] at (P38) {$F$};
  			\node[left]  at (P8) {$G\,\,$};
  			\node[right] at (P10) {$\,\,H$};
  			
  			\spy[black] on (P42) in node at (2.0,-1.25);
  			\spy[black] on (P46) in node at (5.0,-1.25);
  			\spy[black] on (P50) in node at (8.0,-1.25);
  			\spy[black] on (P54) in node at (11.0,-1.25);
  			\end{tikzpicture}
  		
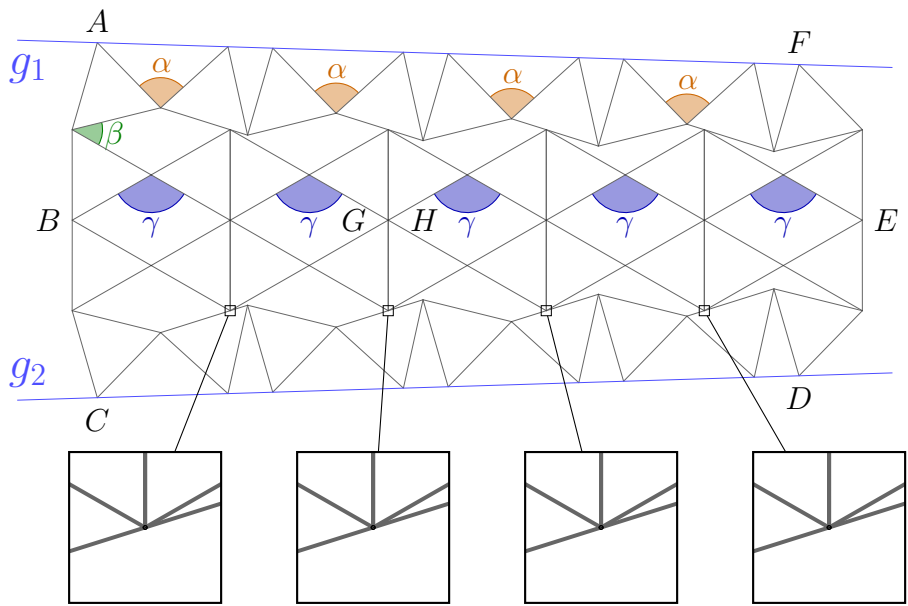
\captionof{figure}{The non-planar subgraph $G_1$.}
  	\end{minipage}
  \end{center}

  For a proof of the non-planarity of $G_1$ we refer the reader to \cite{Harborth}. The proof for its existence is similar to the proof of the main Theorem.
  
  \section{\large{The main Theorem }}
  
  For the main Theorem we use the planar subgraph $G_2$ which also consists of 38 unit triangles (see Figure 3). The graph achieves its planarity through a slight modification of the right side.
  
  \begin{center}
  	\begin{minipage}{\linewidth}
  	  \centering
  	  
  	  \xdef\LstPN{0}
  	  \newif\ifDupe
  	  \pgfplotsset{avoid dupes/.code={\Dupefalse
  	  		\xdef\anker{\DefaultTextposition} 
  	  		\foreach \X in \LstPN
  	  		{\pgfmathtruncatemacro{\itest}{ifthenelse(\X==\punktnummer,1,0)}
  	  			\ifnum\itest=1
  	  			\global\Dupetrue
  	  			\breakforeach
  	  			\fi}
  	  		\ifDupe
  	  		\typeout{\punktnummer\space ist\space ein\space Duplikat!}%
  	  		\xdef\punktnummer{} 
  	  		\else
  	  		\xdef\LstPN{\LstPN,\punktnummer}
  	  		\typeout{\punktnummer\space ist\space neu\space mit\space urprgl.\space Anker=\anker}
  	  		\foreach \X in \LstExcept
  	  		{\ifnum\X=\punktnummer
  	  			\xdef\anker{\AusnahmeTextposition}
  	  			\fi}
  	  		\typeout{\punktnummer\space ist\space neu\space mit\space Anker=\anker}
  	  		\fi}}
    		
  	  \def\DefaultTextposition{east}
  	  \def\AusnahmeTextposition{west}
  	  \def\AusnahmeListe{}
  	  \xdef\BeliebigesVorhandenesKoordinatenpaar{{0.86604121676819312281, 2.42986917293281434738}} 
  	  \colorlet{Kantenfarbe}{gray}
  	  \colorlet{Punktfarbe}{black}
  	  \def\Beschriftung{} 
  	  \pgfplotsset{
  	  	x=12.0mm, y=12.0mm,  
  	  }
  	  
  	  	\xdef\LstExcept{\AusnahmeListe}
  	  	\pgfdeclarelayer{bg}    
  	  	\pgfsetlayers{bg,main}  
  	  	
  	  	\pgfmathsetmacro{\xAlias}{\BeliebigesVorhandenesKoordinatenpaar[0]}
  	  	\pgfmathsetmacro{\yAlias}{\BeliebigesVorhandenesKoordinatenpaar[1]}
  	  	
  	  	\begin{tikzpicture}[SpyStyle]
  	  	
  	  	
  	  	\begin{axis}[hide axis, colormap={kantenfarbe}{color=(Kantenfarbe) color=(Kantenfarbe)}, line width=0.1, 
  	  	]
  	  	\addplot+[mark size=0.01pt, 
  	  	mark options={Punktfarbe}, 
  	  	table/row sep=newline, 
  	  	patch, 
  	  	patch type=polygon,
  	  	vertex count=2, 
  	  	%
  	  	patch table with point meta={
  	  		Startpkt Endpkt colordata  \\
  	  		1 1 \\
  	  		2 1 \\
  	  		3 1 \\
  	  		3 2 \\
  	  		3 38 \\
  	  		3 39 \\
  	  		4 1 \\
  	  		5 1 \\
  	  		5 4 \\
  	  		6 5 \\
  	  		7 5 \\
  	  		7 6 \\
  	  		8 7 \\
  	  		9 7 \\
  	  		9 8 \\
  	  		10 9 \\
  	  		11 9 \\
  	  		11 10 \\
  	  		12 11 \\
  	  		13 11 \\
  	  		13 12 \\
  	  		14 13 \\
  	  		14 49 \\
  	  		15 13 \\
  	  		15 14 \\
  	  		16 2 \\
  	  		17 2 \\
  	  		17 16 \\
  	  		18 16 \\
  	  		19 18 \\
  	  		19 16 \\
  	  		20 19 \\
  	  		21 20 \\
  	  		21 19 \\
  	  		22 21 \\
  	  		23 22 \\
  	  		23 21 \\
  	  		24 23 \\
  	  		25 24 \\
  	  		25 23 \\
  	  		26 25 \\
  	  		27 26 \\
  	  		27 25 \\
  	  		28 27 \\
  	  		29 28 \\
  	  		29 27 \\
  	  		29 15 \\
  	  		30 29 \\
  	  		30 15 \\
  	  		31 30 \\
  	  		32 31 \\
  	  		32 30 \\
  	  		33 31 \\
  	  		34 33 \\
  	  		34 31 \\
  	  		35 34 \\
  	  		35 32 \\
  	  		36 34 \\
  	  		36 35 \\
  	  		37 35 \\
  	  		37 32 \\
  	  		38 38 \\
  	  		39 38 \\
  	  		40 38 \\
  	  		41 38 \\
  	  		41 40 \\
  	  		42 41 \\
  	  		43 41 \\
  	  		43 42 \\
  	  		44 43 \\
  	  		45 43 \\
  	  		45 44 \\
  	  		46 45 \\
  	  		47 45 \\
  	  		47 46 \\
  	  		48 47 \\
  	  		49 47 \\
  	  		49 48 \\
  	  		50 49 \\
  	  		50 14 \\
  	  		51 39 \\
  	  		52 39 \\
  	  		52 51 \\
  	  		53 51 \\
  	  		54 51 \\
  	  		54 53 \\
  	  		55 54 \\
  	  		56 54 \\
  	  		56 55 \\
  	  		57 56 \\
  	  		58 56 \\
  	  		58 57 \\
  	  		59 58 \\
  	  		60 58 \\
  	  		60 59 \\
  	  		61 60 \\
  	  		62 60 \\
  	  		62 61 \\
  	  		63 62 \\
  	  		64 50 \\
  	  		64 62 \\
  	  		64 63 \\
  	  		65 50 \\
  	  		65 64 \\
  	  		66 65 \\
  	  		67 65 \\
  	  		67 66 \\
  	  		68 66 \\
  	  		69 66 \\
  	  		69 68 \\
  	  		70 67 \\
  	  		70 69 \\
  	  		71 69 \\
  	  		71 70 \\
  	  		72 67 \\
  	  		72 70 \\
  	  	},
  	  	%
  	  	visualization depends on={value \thisrowno{0} \as \punktnummer},
  	  	every node near coord/.append style={
  	  		/pgfplots/avoid dupes,
  	  	},
  	  	nodes near coords={\Beschriftung},
  	  	nodes near coords style={
  	  		anchor=\anker,
  	  		text=black, font=\scriptsize, 
  	  		name=p-\punktnummer, 
  	  		path picture={
  	  			\coordinate[] (P\punktnummer) at (p-\punktnummer.\anker);}
  	  	},
  	  	]
  	  	table[header=true, x index=1, y index=2, row sep=\\] {
  	  		Nr x y        \\
  	  		0 0 0         \\
  	  		1 0.86604121676819312281 2.42986917293281434738  \\
  	  		2 0.00000000000000318341 2.92984178304128128900  \\
  	  		3 0.00003162683375346019 1.92984178354140967215  \\
  	  		4 1.73205080670263167342 1.92984178354140767375  \\
  	  		5 1.73208243353638402517 2.92984178304127906856  \\
  	  		6 1.73211406037013393444 1.92984178354140767375  \\
  	  		7 2.59812365030457348425 2.42986917293281212693  \\
  	  		8 3.46413324023901170179 1.92984178354140545331  \\
  	  		9 3.46416486707276405355 2.92984178304127684811  \\
  	  		10 3.46419649390651329668 1.92984178354140545331  \\
  	  		11 4.33020608384095329058 2.42986917293280946240  \\
  	  		12 5.19621567377539150812 1.92984178354140190059  \\
  	  		13 5.19624730060914430396 2.92984178304127329540  \\
  	  		14 5.19627892744289354710 1.92984178354140190059  \\
  	  		15 6.06228851737733354099 2.42986917293280590968  \\
  	  		16 0.98924639929076829681 3.07610032634818164610  \\
  	  		17 0.36795958562110325785 3.85968356708281890022  \\
  	  		18 1.63922970161370340492 3.83604881585347179396  \\
  	  		19 1.97237274793142880469 2.89317251925346363706  \\
  	  		20 2.34033233355252878738 3.82301430329500036009  \\
  	  		21 2.96161914722219288265 3.03943106256036443824  \\
  	  		22 3.61160244954512821280 3.79937955206565369792  \\
  	  		23 3.94474549586285228031 2.85650325546564642920  \\
  	  		24 4.31270508148395315118 3.78634503950718315224  \\
  	  		25 4.93399189515361680236 3.00276179877254678630  \\
  	  		26 5.58397519747655302069 3.76271028827783649007  \\
  	  		27 5.91711824379427753229 2.81983399167782877726  \\
  	  		28 6.28507782941537840316 3.74967577571936594438  \\
  	  		29 6.90636464308504205434 2.96609253498472957844  \\
  	  		30 6.94870963387085360807 1.96698948636794690437  \\
  	  		31 7.43063077127556503854 2.84320408184245065897  \\
  	  		32 7.94849430142083601680 1.98774083649202792934  \\
  	  		33 6.98161160061860286419 3.73672620489399243127  \\
  	  		34 7.97993404335312384035 3.67882715194337217568  \\
  	  		35 8.49779757349839393044 2.82336390659294922401  \\
  	  		36 8.97971871090310358454 3.69957850206745320065  \\
  	  		37 8.94681674415535610478 1.92984178354140767375  \\
  	  		38 0.86604121676819167952 1.42981439415000233240  \\
  	  		39 0.00000000000000000000 0.92984178404153827735  \\
  	  		40 1.73205080670263167342 1.92984178354140656353  \\
  	  		41 1.73208243353638091655 0.92984178404153539077  \\
  	  		42 1.73211406037013393444 1.92984178354140656353  \\
  	  		43 2.59812365030457215198 1.42981439414999944582  \\
  	  		44 3.46413324023901170179 1.92984178354140367695  \\
  	  		45 3.46416486707276094492 0.92984178404153183806  \\
  	  		46 3.46419649390651329668 1.92984178354140367695  \\
  	  		47 4.33020608384095151422 1.42981439414999700332  \\
  	  		48 5.19621567377539150812 1.92984178354140234468  \\
  	  		49 5.19624730060914075125 0.92984178404153117192  \\
  	  		50 6.06228851737733087646 1.42981439414999456083  \\
  	  		51 0.98924639929076485512 0.78358324073463470061  \\
  	  		52 0.36795958562109704060 0.00000000000000000000  \\
  	  		53 1.63922970161369718767 0.02363475122934224557  \\
  	  		54 1.97237274793142547402 0.96651104782934993409  \\
  	  		55 2.34033233355252345831 0.03666926378781244084  \\
  	  		56 2.96161914722218977403 0.82025250452244613530  \\
  	  		57 3.61160244954512243964 0.06030401501715468987  \\
  	  		58 3.94474549586284961578 1.00318031161716070265  \\
  	  		59 4.31270508148394782211 0.07333852757562275837  \\
  	  		60 4.93399189515361413783 0.85692176831025679284  \\
  	  		61 5.58397519747654769162 0.09697327880496642294  \\
  	  		62 5.91711824379427575593 1.03984957540497280348  \\
  	  		63 6.28507782941537396226 0.11000779136343377673  \\
  	  		64 6.90636464308503938980 0.89359103209806989288  \\
  	  		65 6.94870963387085360807 1.89269408071485112366  \\
  	  		66 7.43063077127556237400 1.01647948524034648088  \\
  	  		67 7.94849430142083690498 1.87194273059076876642  \\
  	  		68 6.98161160061859664694 0.12295736218880580493  \\
  	  		69 7.97993404335311939946 0.18085641513942235514  \\
  	  		70 8.49779757349839037772 1.03631966048984347495  \\
  	  		71 8.97971871090310003183 0.16010506501534016444  \\
  	  		72 8.94681674415535610478 1.92984178354138413702  \\
  	  	};
  	  	
  	  	\addplot[no marks, 
  	  	nodes near coords={},
  	  	visualization depends on={value \thisrowno{0} \as \PunktI},
  	  	visualization depends on={value \thisrowno{1} \as \Scheitel},
  	  	visualization depends on={value \thisrowno{2} \as \PunktII},
  	  	visualization depends on={value \thisrowno{3} \as \Winkelradius},
  	  	visualization depends on={value \thisrowno{4} \as \Winkelfarbe},
  	  	visualization depends on={value \thisrowno{5} \as \Winkelname},
  	  	visualization depends on={value \thisrowno{6} \as \WinkelExzentrizitaet},
  	  	nodes near coords style={anchor=center,
  	  		path picture={
  	  			\begin{pgfonlayer}{bg}    
  	  			\draw pic [angle radius=\Winkelradius cm,%
  	  			fill=\Winkelfarbe!40, draw=\Winkelfarbe,
  	  			"$\Winkelname$", angle eccentricity =\WinkelExzentrizitaet,
  	  			text=\Winkelfarbe%
  	  			] {angle = P\PunktI--P\Scheitel--P\PunktII};
  	  			\end{pgfonlayer}
  	  	}},%
  	  	]
  	  	table[header=true, x expr=0.86604121676819312281, y expr=2.42986917293281434738]{
  	  		Punkt1 Scheitel Punkt2 Winkelradius[cm] Winkelfarbe Winkelname WinkelExz
  	  		3 1 4 0.5 blue!70!black {\gamma} 1.4 \\
  	  		6 7 8 0.5 blue!70!black \gamma{} 1.4 \\
  	  		10 11 12 0.5 blue {\gamma} 1.4 \\
  	  		1 2 16 0.4 green!50!black {\beta} 1.4 \\
  	  		18 16 17 0.4 orange!80!black {\alpha} 1.4 \\
  	  		22 21 20 0.4 orange!80!black {\alpha} 1.4 \\
  	  		26 25 24 0.4 orange!80!black {\alpha} 1.4 \\
  	  		14 15 30 0.5 violet {\delta} 1.4 \\
  	  	};
  	  	
  	  	\end{axis}
  	  	
  	  	\draw[blue!70!white, line width=0.1, shorten >= -15, shorten <=-30] (P17) -- (P36);
  	  	\draw[blue!70!white, line width=0.1, shorten >= -15, shorten <=-30] (P52) -- (P71);
  	  	\node[below, blue!70!white, font=\Large] at (P17) {$g_{1}\quad\quad\quad$};
  	  	\node[above, blue!70!white, font=\Large] at (P52) {$g_{2}\quad\quad\quad$};
  	  	
  	  	\node[above] at (P17) {$A$};
  	  	\node[left]  at (P3) {$B$};
  	  	\node[below] at (P52) {$C$};
  	  	\node[below] at (P71) {$D$};
  	  	\node[right] at (P37) {$E$};
  	  	\node[above] at (P36) {$F$};
  	  	\node[left]  at (P8) {$G\,\,$};
  	  	\node[right] at (P10) {$\,\,H$};
  	  	
  	  	\spy[black] on (P41) in node at (2.0,-1.25);
  	  	\spy[black] on (P45) in node at (5.0,-1.25);
  	  	\spy[black] on (P49) in node at (8.0,-1.25);
  	  	\spy[black] on (P69) in node at (11.0,-1.25);
  	  	\end{tikzpicture}
  	  
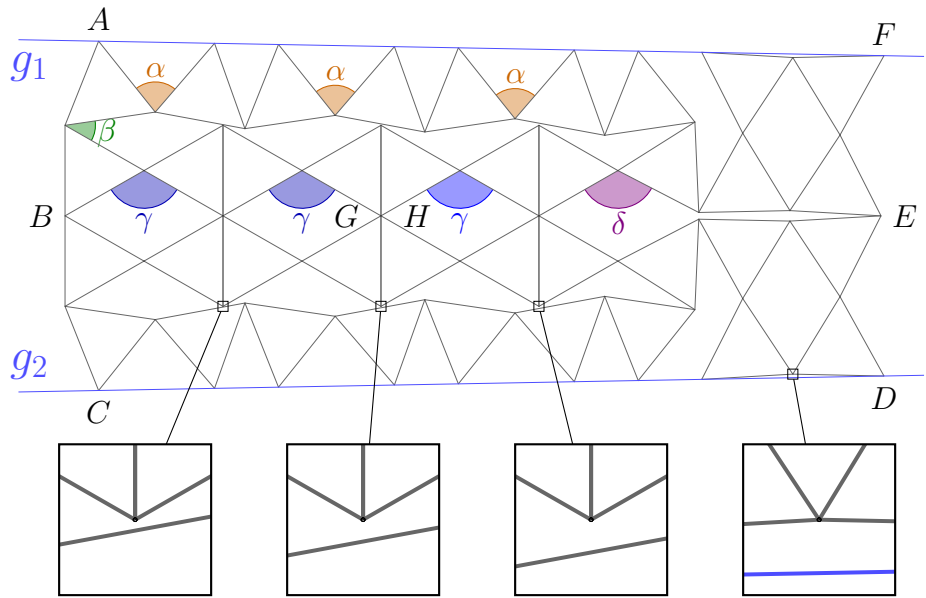
\captionof{figure}{The planar subgraph $G_2$.}
  	\end{minipage}
  \end{center}
  \quad\\
  
  \textbf{Theorem.} \textit{The number of unit triangles in 4-regular planar unit-distance graphs consisting only of unit triangles without additional triangles is $\leq6422$.}
  \\
  
  \textit{Proof.} We arrange 38 unit triangles in the plane in such a way that we get the graph $G_2$ which only contains vertices of degree 2 and 4. The vertices of degree 2 lie exactly on the two intersecting lines $g_1$ and $g_2$, where the point of intersection is right of $E$ (see Figure 3). $G_2$ is flexible and mirror-symmetric about $\overline{BE}$. Let $n$ be an integer and let $\omega=\frac{360^\circ}{n}$ be the angle between $g_1$ and $g_2$. Let us imagine $g_1$ and $g_2$ as rails on which the vertices with the exception of $D$ and $F$ are movably mounted. Then, because of its flexibility, $G_2$ can be pushed back and forth slightly, whereby $D$ and $F$ move up and down. With such a movement no angle in the graph remains constant. If one of the angles $\alpha$, $\beta$, or $\gamma$ is changed continuously close to its value given in Table 1, then the path of $D$ and $F$ sometime cuts $g_1$ and $g_2$. This intersection can not be specified exactly, but it must exist. During this movement of $D$ and $F$ across $g_1$ and $g_2$ it must be ensured that each vertex of degree 4 do not intersect with other edges nor with $g_1$ and $g_2$. The smallest integer $n$ for which this applies is $n=169$. Because of the mirror symmetry we can connect 169 copies of $G_2$ together in such a way that the vertices $A$ and $C$, as well as $F$ and $D$ coincide to build a rigid planar 4-regular ring graph consisting of $169\cdot38=6422$ unit triangles. The vertex of $\omega$ is the center of this ring.\hfill$\square$
  \\
  
  Because $G_1$ is also mirror-symmetric about $\overline{BE}$, the proof for the existence of a rigid ring graph consisting of 100 copies of $G_1$  is similar to the proof before. Table 1 gives the approximate values for the angles and distances in both subgraphs.
  \\ \\
  
  \begin{center}
  	\begin{tabular}[t]{r|r|r}
  		& $G_1$ \quad\quad\quad\quad\quad & $G_2$ \quad\quad\quad\quad\quad \\
  		\hline & & \\
  		$\measuredangle(\overrightarrow{FA},\overrightarrow{DC})=$ & $360^\circ/100=3{,}6^\circ$ &  $360^\circ/169\approx2{,}130^\circ$ \\
  		$\alpha\approx$        &  $91{,}58566772584003^\circ$ &  $78{,}95050838942406^\circ$ \\
  		$\beta\approx$         &  $43{,}94364026236698^\circ$ & $38{,}40835335322197^\circ$ \\
  		$\gamma\approx$        & $119{,}90161279889431^\circ$ &  $119{,}99637583181277^\circ$ \\
  		$\delta\approx$        & -- \quad\quad\quad\quad\quad & $122{,}42510282308054^\circ$ \\
  		$\overline{GH}\approx$ &   $0{,}00171718039014 \;\,$  &   $0{,}00006325366750 \;\,$
  	\end{tabular}
  	\captionof{table}{}
  \end{center}
  
  \newpage
  
  There exists also a 4-regular planar unit triangle graph with a rotational symmetry of order 100. This rigid ring graph consists of 9200 unit triangles and is build from the subgraph $G_3$ (see Figure 4). As well as $G_2$, these graph contains very small distances between vertices and edges.
  \\
  
  \begin{center}
  	\begin{minipage}{\linewidth}
  	  \centering	
  	  \begin{tikzpicture}
  	[y=0.22pt, x=0.22pt, yscale=-1.0, xscale=1.0]
  	
  	\draw[black, line width=0.01pt]
  	(38.4317,915.0584) -- (76.3458,936.9299)
  	
  	(38.4475,958.8287) -- (76.3458,936.9299)
  	
  	(38.4475,958.8287) -- (38.4317,915.0584)
  	
  	(38.4475,958.8287) -- (76.3458,980.7276)
  	
  	(38.4475,958.8287) -- (38.4317,1002.5991)
  	
  	(114.2442,958.8287) -- (76.3458,936.9299)
  	
  	(114.2442,958.8287) -- (76.3458,980.7276)
  	
  	(114.2442,958.8287) -- (114.2600,1002.5991)
  	
  	(114.2600,915.0584) -- (76.3458,936.9299)
  	
  	(114.2600,915.0584) -- (114.2442,958.8287)
  	
  	(114.2757,958.8287) -- (114.2599,915.0584)
  	
  	(114.2757,958.8287) -- (114.2599,1002.5991)
  	
  	(114.2757,958.8287) -- (152.1741,980.7276)
  	
  	(152.1741,936.9299) -- (114.2600,915.0584)
  	
  	(152.1741,936.9299) -- (114.2757,958.8287)
  	
  	(190.0724,958.8287) -- (152.1741,936.9299)
  	
  	(190.0724,958.8287) -- (152.1741,980.7276)
  	
  	(190.0724,958.8287) -- (190.0882,1002.5991)
  	
  	(190.0882,915.0584) -- (152.1741,936.9299)
  	
  	(190.0882,915.0584) -- (190.0724,958.8287)
  	
  	(190.1040,958.8287) -- (190.0882,915.0584)
  	
  	(190.1040,958.8287) -- (190.0882,1002.5991)
  	
  	(190.1040,958.8287) -- (228.0023,980.7276)
  	
  	(228.0023,936.9299) -- (190.0882,915.0584)
  	
  	(228.0023,936.9299) -- (190.1040,958.8287)
  	
  	(265.9007,958.8287) -- (228.0023,936.9299)
  	
  	(265.9007,958.8287) -- (228.0023,980.7276)
  	
  	(265.9007,958.8287) -- (265.9165,1002.5991)
  	
  	(265.9165,915.0584) -- (228.0023,936.9299)
  	
  	(265.9165,915.0584) -- (265.9006,958.8287)
  	
  	(265.9322,958.8287) -- (265.9164,915.0584)
  	
  	(265.9322,958.8287) -- (265.9164,1002.5991)
  	
  	(265.9322,958.8287) -- (303.8306,980.7276)
  	
  	(303.8306,936.9299) -- (265.9165,915.0584)
  	
  	(303.8306,936.9299) -- (265.9322,958.8287)
  	
  	(341.7289,958.8287) -- (303.8306,936.9299)
  	
  	(341.7289,958.8287) -- (303.8306,980.7276)
  	
  	(341.7289,958.8287) -- (341.7447,1002.5991)
  	
  	(341.7447,915.0584) -- (303.8306,936.9299)
  	
  	(341.7447,915.0584) -- (341.7289,958.8287)
  	
  	(341.7605,958.8287) -- (341.7447,915.0584)
  	
  	(341.7605,958.8287) -- (341.7447,1002.5991)
  	
  	(379.6588,936.9299) -- (341.7447,915.0584)
  	
  	(379.6588,936.9299) -- (341.7605,958.8287)
  	
  	(379.6588,936.9299) -- (417.5371,914.9964)
  	
  	(80.8887,904.4168) -- (38.4317,915.0584)
  	
  	(50.4443,872.9687) -- (38.4317,915.0584)
  	
  	(50.4443,872.9687) -- (80.8887,904.4168)
  	
  	(113.2477,874.9424) -- (80.8887,904.4168)
  	
  	(122.5938,917.7033) -- (113.2477,874.9424)
  	
  	(122.5938,917.7033) -- (80.8887,904.4168)
  	
  	(134.6064,875.6136) -- (122.5938,917.7033)
  	
  	(165.0508,907.0617) -- (134.6064,875.6136)
  	
  	(165.0508,907.0617) -- (122.5938,917.7033)
  	
  	(197.4098,877.5873) -- (165.0508,907.0617)
  	
  	(206.7559,920.3482) -- (197.4098,877.5873)
  	
  	(206.7559,920.3482) -- (165.0508,907.0617)
  	
  	(218.7685,878.2585) -- (206.7559,920.3482)
  	
  	(249.2129,909.7066) -- (218.7685,878.2585)
  	
  	(249.2129,909.7066) -- (206.7559,920.3482)
  	
  	(281.5719,880.2322) -- (249.2129,909.7066)
  	
  	(290.9180,922.9931) -- (281.5719,880.2322)
  	
  	(290.9180,922.9931) -- (249.2129,909.7066)
  	
  	(302.9306,880.9034) -- (290.9180,922.9931)
  	
  	(333.3750,912.3515) -- (302.9306,880.9034)
  	
  	(333.3750,912.3515) -- (290.9180,922.9931)
  	
  	(365.7340,882.8771) -- (333.3750,912.3515)
  	
  	(375.0801,925.6380) -- (365.7340,882.8771)
  	
  	(375.0801,925.6380) -- (333.3750,912.3515)
  	
  	(387.0927,883.5483) -- (375.0801,925.6380)
  	
  	(417.5371,914.9964) -- (387.0927,883.5483)
  	
  	(417.5371,914.9964) -- (375.0801,925.6380)
  	
  	(417.5929,958.7667) -- (417.5370,914.9964)
  	
  	(417.5929,958.7667) -- (379.6588,936.9299)
  	
  	(461.3633,958.8287) -- (417.5930,958.7667)
  	
  	(461.3633,958.8287) -- (417.5930,958.8908)
  	
  	(461.3633,958.8287) -- (439.5319,996.7660)
  	
  	(439.5319,920.8915) -- (417.5930,958.7667)
  	
  	(439.5319,920.8915) -- (461.3633,958.8287)
  	
  	(465.9965,886.0280) -- (439.5319,920.8915)
  	
  	(482.9569,926.3788) -- (465.9965,886.0280)
  	
  	(482.9569,926.3788) -- (439.5318,920.8915)
  	
  	(512.3311,958.8287) -- (482.9569,926.3788)
  	
  	(512.3311,958.8287) -- (482.9569,991.2787)
  	
  	(512.3311,958.8287) -- (525.7465,1000.4925)
  	
  	(525.7465,917.1650) -- (482.9569,926.3788)
  	
  	(525.7465,917.1650) -- (512.3311,958.8287)
  	
  	(539.1619,958.8287) -- (525.7465,917.1650)
  	
  	(539.1619,958.8287) -- (525.7465,1000.4925)
  	
  	(539.1619,958.8287) -- (568.5361,991.2787)
  	
  	(568.5361,926.3788) -- (525.7465,917.1650)
  	
  	(568.5361,926.3788) -- (539.1619,958.8287)
  	
  	(592.8999,890.0161) -- (568.5361,926.3788)
  	
  	(612.2090,929.2971) -- (592.8999,890.0161)
  	
  	(612.2090,929.2971) -- (568.5361,926.3788)
  	
  	(644.5158,958.8287) -- (612.2090,929.2971)
  	
  	(644.5158,958.8287) -- (612.2090,988.3604)
  	
  	(644.5158,958.8287) -- (653.9375,1001.5730)
  	
  	(653.9375,916.0844) -- (612.2090,929.2971)
  	
  	(653.9375,916.0844) -- (644.5158,958.8287)
  	
  	(691.1906,893.1050) -- (653.9375,916.0844)
  	
  	(692.4648,936.8568) -- (691.1905,893.1050)
  	
  	(692.4648,936.8568) -- (653.9375,916.0844)
  	
  	(730.3208,958.8287) -- (692.4648,936.8568)
  	
  	(730.3208,958.8287) -- (692.4648,980.8007)
  	
  	(730.3208,958.8287) -- (730.4211,1002.5990)
  	
  	(730.4211,915.0585) -- (692.4648,936.8568)
  	
  	(730.4211,915.0585) -- (730.3208,958.8287)
  	
  	(730.5214,958.8287) -- (730.4211,915.0585)
  	
  	(730.5214,958.8287) -- (730.4211,1002.5990)
  	
  	(730.5214,958.8287) -- (768.3774,980.8007)
  	
  	(768.3774,936.8568) -- (730.4211,915.0585)
  	
  	(768.3774,936.8568) -- (730.5214,958.8287)
  	
  	(806.2334,958.8287) -- (768.3774,936.8568)
  	
  	(806.2334,958.8287) -- (768.3774,980.8007)
  	
  	(806.2334,958.8287) -- (806.3337,1002.5990)
  	
  	(806.3337,915.0585) -- (768.3774,936.8568)
  	
  	(806.3337,915.0585) -- (806.2334,958.8287)
  	
  	(806.4340,958.8287) -- (806.3337,915.0585)
  	
  	(806.4340,958.8287) -- (806.3337,1002.5990)
  	
  	(806.4340,958.8287) -- (844.2900,980.8007)
  	
  	(844.2900,936.8568) -- (806.3337,915.0585)
  	
  	(844.2900,936.8568) -- (806.4340,958.8287)
  	
  	(865.5271,898.5837) -- (844.2900,936.8568)
  	
  	(888.0540,936.1122) -- (865.5271,898.5837)
  	
  	(888.0540,936.1122) -- (844.2900,936.8568)
  	
  	(912.8928,900.0723) -- (888.0540,936.1122)
  	
  	(931.6849,939.6033) -- (912.8928,900.0723)
  	
  	(931.6849,939.6033) -- (888.0540,936.1122)
  	
  	(952.9221,901.3302) -- (931.6849,939.6033)
  	
  	(975.4489,938.8587) -- (952.9221,901.3302)
  	
  	(975.4489,938.8587) -- (931.6849,939.6033)
  	
  	(1000.2877,902.8188) -- (975.4489,938.8587)
  	
  	(1019.0798,942.3498) -- (1000.2877,902.8188)
  	
  	(1019.0798,942.3498) -- (975.4489,938.8587)
  	
  	(1040.3170,904.0767) -- (1019.0798,942.3498)
  	
  	(1062.8438,941.6052) -- (1040.3170,904.0767)
  	
  	(1062.8438,941.6052) -- (1019.0798,942.3498)
  	
  	(1087.6826,905.5652) -- (1062.8438,941.6052)
  	
  	(1106.4747,945.0963) -- (1087.6826,905.5652)
  	
  	(1106.4747,945.0963) -- (1062.8438,941.6052)
  	
  	(1127.7119,906.8232) -- (1106.4747,945.0963)
  	
  	(1150.2387,944.3517) -- (1127.7119,906.8232)
  	
  	(1150.2387,944.3517) -- (1106.4747,945.0963)
  	
  	(1175.0776,908.3117) -- (1150.2387,944.3517)
  	
  	(1193.8696,947.8428) -- (1175.0776,908.3117)
  	
  	(1193.8696,947.8428) -- (1150.2387,944.3517)
  	
  	(1215.1068,909.5697) -- (1193.8696,947.8428)
  	
  	(1237.6337,947.0982) -- (1215.1068,909.5697)
  	
  	(1237.6337,947.0982) -- (1193.8696,947.8428)
  	
  	(1262.4725,911.0582) -- (1237.6337,947.0982)
  	
  	(1281.2646,950.5893) -- (1262.4725,911.0582)
  	
  	(1281.2646,950.5893) -- (1237.6337,947.0982)
  	
  	(1302.5017,912.3162) -- (1281.2646,950.5893)
  	
  	(1325.0286,949.8447) -- (1302.5017,912.3162)
  	
  	(1325.0286,949.8447) -- (1281.2646,950.5893)
  	
  	(1349.8674,913.8047) -- (1325.0286,949.8447)
  	
  	(1368.6595,953.3357) -- (1349.8674,913.8047)
  	
  	(1368.6595,953.3357) -- (1325.0286,949.8447)
  	
  	(1389.8967,915.0627) -- (1368.6595,953.3357)
  	
  	(1412.4235,952.5911) -- (1389.8967,915.0627)
  	
  	(1412.4235,952.5911) -- (1368.6595,953.3357)
  	
  	(1437.2623,916.5512) -- (1412.4235,952.5912)
  	
  	(1456.0544,956.0823) -- (1437.2623,916.5512)
  	
  	(1456.0544,956.0823) -- (1412.4235,952.5912)
  	
  	(1477.2916,917.8092) -- (1456.0544,956.0823)
  	
  	(1499.8184,955.3377) -- (1477.2916,917.8092)
  	
  	(1499.8184,955.3377) -- (1456.0544,956.0823)
  	
  	(1524.6573,919.2977) -- (1499.8184,955.3377)
  	
  	(1543.4493,958.8287) -- (1524.6573,919.2977)
  	
  	(1543.4493,958.8287) -- (1499.8184,955.3377)
  	
  	(1543.4493,958.8287) -- (1499.8184,962.3198)
  	
  	(1543.4493,958.8287) -- (1524.6573,998.3598)
  	
  	(38.4317,1002.5991) -- (76.3458,980.7276)
  	
  	(114.2600,1002.5991) -- (76.3458,980.7276)
  	
  	(152.1741,980.7276) -- (114.2600,1002.5991)
  	
  	(190.0882,1002.5991) -- (152.1741,980.7276)
  	
  	(228.0023,980.7276) -- (190.0882,1002.5991)
  	
  	(265.9165,1002.5991) -- (228.0023,980.7276)
  	
  	(303.8306,980.7276) -- (265.9165,1002.5991)
  	
  	(341.7447,1002.5991) -- (303.8306,980.7276)
  	
  	(379.6588,980.7276) -- (341.7447,1002.5991)
  	
  	(379.6588,980.7276) -- (341.7605,958.8287)
  	
  	(379.6588,980.7276) -- (417.5371,1002.6611)
  	
  	(80.8887,1013.2407) -- (38.4317,1002.5991)
  	
  	(50.4443,1044.6888) -- (38.4317,1002.5991)
  	
  	(50.4443,1044.6888) -- (80.8887,1013.2407)
  	
  	(113.2477,1042.7151) -- (80.8887,1013.2407)
  	
  	(122.5938,999.9542) -- (80.8887,1013.2407)
  	
  	(122.5938,999.9542) -- (113.2477,1042.7151)
  	
  	(134.6064,1042.0439) -- (122.5938,999.9542)
  	
  	(165.0508,1010.5958) -- (122.5938,999.9542)
  	
  	(165.0508,1010.5958) -- (134.6064,1042.0439)
  	
  	(197.4098,1040.0702) -- (165.0508,1010.5958)
  	
  	(206.7559,997.3093) -- (165.0508,1010.5958)
  	
  	(206.7559,997.3093) -- (197.4098,1040.0702)
  	
  	(218.7685,1039.3990) -- (206.7559,997.3093)
  	
  	(249.2129,1007.9509) -- (206.7559,997.3093)
  	
  	(249.2129,1007.9509) -- (218.7685,1039.3990)
  	
  	(281.5719,1037.4253) -- (249.2129,1007.9509)
  	
  	(290.9180,994.6644) -- (249.2129,1007.9509)
  	
  	(290.9180,994.6644) -- (281.5719,1037.4253)
  	
  	(302.9306,1036.7541) -- (290.9180,994.6644)
  	
  	(333.3750,1005.3060) -- (290.9180,994.6644)
  	
  	(333.3750,1005.3060) -- (302.9306,1036.7541)
  	
  	(365.7340,1034.7804) -- (333.3750,1005.3060)
  	
  	(375.0801,992.0195) -- (333.3750,1005.3060)
  	
  	(375.0801,992.0195) -- (365.7340,1034.7804)
  	
  	(387.0927,1034.1092) -- (375.0801,992.0195)
  	
  	(417.5371,1002.6611) -- (375.0801,992.0195)
  	
  	(417.5371,1002.6611) -- (387.0927,1034.1092)
  	
  	(417.5929,958.8908) -- (379.6588,980.7276)
  	
  	(417.5929,958.8908) -- (417.5370,1002.6611)
  	
  	(439.5319,996.7660) -- (417.5929,958.8908)
  	
  	(465.9965,1031.6296) -- (439.5319,996.7660)
  	
  	(482.9569,991.2787) -- (439.5318,996.7660)
  	
  	(482.9569,991.2787) -- (465.9965,1031.6296)
  	
  	(525.7465,1000.4925) -- (482.9569,991.2787)
  	
  	(568.5361,991.2787) -- (525.7465,1000.4925)
  	
  	(592.8999,1027.6414) -- (568.5361,991.2787)
  	
  	(612.2090,988.3604) -- (568.5361,991.2787)
  	
  	(612.2090,988.3604) -- (592.8999,1027.6414)
  	
  	(653.9375,1001.5730) -- (612.2090,988.3604)
  	
  	(691.1906,1024.5525) -- (653.9375,1001.5730)
  	
  	(692.4648,980.8007) -- (653.9375,1001.5730)
  	
  	(692.4648,980.8007) -- (691.1905,1024.5525)
  	
  	(730.4211,1002.5990) -- (692.4648,980.8007)
  	
  	(768.3774,980.8007) -- (730.4211,1002.5990)
  	
  	(806.3337,1002.5990) -- (768.3774,980.8007)
  	
  	(844.2900,980.8007) -- (806.3337,1002.5990)
  	
  	(865.5271,1019.0738) -- (844.2900,980.8007)
  	
  	(888.0540,981.5453) -- (844.2900,980.8007)
  	
  	(888.0540,981.5453) -- (865.5271,1019.0738)
  	
  	(912.8928,1017.5853) -- (888.0540,981.5453)
  	
  	(931.6849,978.0542) -- (888.0540,981.5453)
  	
  	(931.6849,978.0542) -- (912.8928,1017.5853)
  	
  	(952.9221,1016.3273) -- (931.6849,978.0542)
  	
  	(975.4489,978.7988) -- (931.6849,978.0542)
  	
  	(975.4489,978.7988) -- (952.9221,1016.3273)
  	
  	(1000.2877,1014.8387) -- (975.4489,978.7988)
  	
  	(1019.0798,975.3077) -- (975.4489,978.7988)
  	
  	(1019.0798,975.3077) -- (1000.2877,1014.8387)
  	
  	(1040.3170,1013.5808) -- (1019.0798,975.3077)
  	
  	(1062.8438,976.0523) -- (1019.0798,975.3077)
  	
  	(1062.8438,976.0523) -- (1040.3170,1013.5808)
  	
  	(1087.6826,1012.0922) -- (1062.8438,976.0523)
  	
  	(1106.4747,972.5612) -- (1062.8438,976.0523)
  	
  	(1106.4747,972.5612) -- (1087.6826,1012.0922)
  	
  	(1127.7119,1010.8343) -- (1106.4747,972.5612)
  	
  	(1150.2387,973.3058) -- (1106.4747,972.5612)
  	
  	(1150.2387,973.3058) -- (1127.7119,1010.8343)
  	
  	(1175.0776,1009.3457) -- (1150.2387,973.3058)
  	
  	(1193.8696,969.8147) -- (1150.2387,973.3058)
  	
  	(1193.8696,969.8147) -- (1175.0776,1009.3457)
  	
  	(1215.1068,1008.0878) -- (1193.8696,969.8147)
  	
  	(1237.6337,970.5593) -- (1193.8696,969.8147)
  	
  	(1237.6337,970.5593) -- (1215.1068,1008.0878)
  	
  	(1262.4725,1006.5992) -- (1237.6337,970.5593)
  	
  	(1281.2646,967.0682) -- (1237.6337,970.5593)
  	
  	(1281.2646,967.0682) -- (1262.4725,1006.5992)
  	
  	(1302.5017,1005.3413) -- (1281.2646,967.0682)
  	
  	(1325.0286,967.8128) -- (1281.2646,967.0682)
  	
  	(1325.0286,967.8128) -- (1302.5017,1005.3413)
  	
  	(1349.8674,1003.8528) -- (1325.0286,967.8128)
  	
  	(1368.6595,964.3217) -- (1325.0286,967.8128)
  	
  	(1368.6595,964.3217) -- (1349.8674,1003.8528)
  	
  	(1389.8967,1002.5948) -- (1368.6595,964.3217)
  	
  	(1412.4235,965.0663) -- (1368.6595,964.3217)
  	
  	(1412.4235,965.0663) -- (1389.8967,1002.5948)
  	
  	(1437.2623,1001.1063) -- (1412.4235,965.0663)
  	
  	(1456.0544,961.5752) -- (1412.4235,965.0663)
  	
  	(1456.0544,961.5752) -- (1437.2623,1001.1063)
  	
  	(1477.2916,999.8483) -- (1456.0544,961.5752)
  	
  	(1499.8184,962.3198) -- (1456.0544,961.5752)
  	
  	(1499.8184,962.3198) -- (1477.2916,999.8483)
  	
  	(1524.6573,998.3598) -- (1499.8184,962.3198);
  	
  	\draw[blue!70!white, line width=0.1pt]
  	(6.6956,871.5939) -- (1568.4060,920.6726)
  	(6.6956,1046.0637) -- (1568.4060,996.9849);
  	
  	\end{tikzpicture}
  	
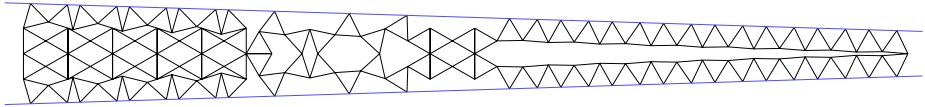
\captionof{figure}{The planar subgraph $G_3$.}
  	\end{minipage}
  \end{center}
  \quad\\
  
  \section{\large{Other forms besides rings}}
  
  In $G_1$ the vertices on $\overline{BE}$ are mirror-symmetric about $g_3$, the perpendicular bisector of $\overline{BE}$. Let $G_4$ be the graph we get by reflecting each vertex beneath $\overline{BE}$ over $g_3$, by which $g_1$ and $g_2$ become parallel (see Figure 5). Let $G_5$ be the graph we get by reflecting each vertex of $G_4$ over $\overline{BE}$. Now we can use $G_4$ and $G_5$ as a kind of adapters to change the direction of the ring segments build with $G_1$, which makes infinitely many other forms besides ring graphs possible (see Figure 6). But, as well as $G_1$, all these graphs are not planar. It remains an open question, if also planar unit triangle graphs in other forms besides rings exist.
  
  \begin{center}
  	\begin{minipage}{\linewidth}
  		\centering
  		
  		\xdef\LstPN{0}
  		\newif\ifDupe
  		\pgfplotsset{avoid dupes/.code={\Dupefalse
  				\xdef\anker{\DefaultTextposition} 
  				\foreach \X in \LstPN
  				{\pgfmathtruncatemacro{\itest}{ifthenelse(\X==\punktnummer,1,0)}
  					\ifnum\itest=1
  					\global\Dupetrue
  					\breakforeach
  					\fi}
  				\ifDupe
  				\typeout{\punktnummer\space ist\space ein\space Duplikat!}%
  				\xdef\punktnummer{} 
  				\else
  				\xdef\LstPN{\LstPN,\punktnummer}
  				\typeout{\punktnummer\space ist\space neu\space mit\space urprgl.\space Anker=\anker}
  				\foreach \X in \LstExcept
  				{\ifnum\X=\punktnummer
  					\xdef\anker{\AusnahmeTextposition}
  					\fi}
  				\typeout{\punktnummer\space ist\space neu\space mit\space Anker=\anker}
  				\fi}}
  			
  		\def\DefaultTextposition{east}
  		\def\AusnahmeTextposition{west}
  		\def\AusnahmeListe{}
  		\xdef\BeliebigesVorhandenesKoordinatenpaar{{0.86645437967485894326, 2.46197022175259894183}} 
  		\colorlet{Kantenfarbe}{gray}
  		\colorlet{Punktfarbe}{black}
  		\def\Beschriftung{} 
  		\pgfplotsset{
  			x=12.0mm, y=12.0mm,  
  		}
  			
  			\xdef\LstExcept{\AusnahmeListe}
  			\pgfdeclarelayer{bg}    
  			\pgfsetlayers{bg,main}  
  			
  			\pgfmathsetmacro{\xAlias}{\BeliebigesVorhandenesKoordinatenpaar[0]}
  			\pgfmathsetmacro{\yAlias}{\BeliebigesVorhandenesKoordinatenpaar[1]}
  			
  			\begin{tikzpicture}[SpyStyle]
  			
  			
  			\begin{axis}[hide axis, colormap={kantenfarbe}{color=(Kantenfarbe) color=(Kantenfarbe)}, line width=0.1, 
  			]
  			\addplot+[mark size=0.01pt, 
  			mark options={Punktfarbe}, 
  			table/row sep=newline, 
  			patch, 
  			patch type=polygon,
  			vertex count=2, 
  			%
  			patch table with point meta={
  				Startpkt Endpkt colordata  \\
  				1 1 \\
  				2 1 \\
  				3 1 \\
  				3 2 \\
  				3 56 \\
  				4 1 \\
  				5 1 \\
  				5 4 \\
  				6 5 \\
  				7 5 \\
  				7 6 \\
  				8 7 \\
  				9 7 \\
  				9 8 \\
  				10 9 \\
  				11 9 \\
  				11 10 \\
  				12 11 \\
  				13 11 \\
  				13 12 \\
  				14 13 \\
  				15 13 \\
  				15 14 \\
  				16 15 \\
  				17 15 \\
  				17 16 \\
  				18 17 \\
  				19 17 \\
  				19 18 \\
  				20 19 \\
  				20 39 \\
  				20 40 \\
  				21 19 \\
  				21 20 \\
  				22 2 \\
  				23 2 \\
  				23 22 \\
  				24 22 \\
  				25 24 \\
  				25 22 \\
  				26 25 \\
  				27 26 \\
  				27 25 \\
  				28 27 \\
  				29 28 \\
  				29 27 \\
  				30 29 \\
  				31 30 \\
  				31 29 \\
  				32 31 \\
  				33 32 \\
  				33 31 \\
  				34 33 \\
  				35 34 \\
  				35 33 \\
  				36 35 \\
  				37 36 \\
  				37 35 \\
  				37 21 \\
  				38 37 \\
  				38 21 \\
  				39 39 \\
  				40 39 \\
  				41 39 \\
  				42 39 \\
  				42 41 \\
  				43 42 \\
  				44 42 \\
  				44 43 \\
  				45 44 \\
  				46 44 \\
  				46 45 \\
  				47 46 \\
  				48 46 \\
  				48 47 \\
  				49 48 \\
  				50 48 \\
  				50 49 \\
  				51 50 \\
  				52 50 \\
  				52 51 \\
  				53 52 \\
  				54 52 \\
  				54 53 \\
  				55 54 \\
  				56 54 \\
  				56 55 \\
  				57 56 \\
  				57 3 \\
  				58 40 \\
  				59 40 \\
  				59 58 \\
  				60 58 \\
  				61 58 \\
  				61 60 \\
  				62 61 \\
  				63 61 \\
  				63 62 \\
  				64 63 \\
  				65 63 \\
  				65 64 \\
  				66 65 \\
  				67 65 \\
  				67 66 \\
  				68 67 \\
  				69 67 \\
  				69 68 \\
  				70 69 \\
  				71 69 \\
  				71 70 \\
  				72 71 \\
  				73 57 \\
  				73 71 \\
  				73 72 \\
  				74 57 \\
  				74 73 \\
  			},
  			%
  			visualization depends on={value \thisrowno{0} \as \punktnummer},
  			every node near coord/.append style={
  				/pgfplots/avoid dupes,
  			},
  			nodes near coords={\Beschriftung},
  			nodes near coords style={
  				anchor=\anker,
  				text=black, font=\scriptsize, 
  				name=p-\punktnummer, 
  				path picture={
  					\coordinate[] (P\punktnummer) at (p-\punktnummer.\anker);}
  			},
  			]
  			table[header=true, x index=1, y index=2, row sep=\\] {
  				Nr x y        \\
  				0 0 0         \\
  				1 0.86645437967485894326 2.46197022175259894183  \\
  				2 0.00000000000000205166 2.96122647653791348787  \\
  				3 0.00085859019507177775 1.96122684512654288369  \\
  				4 1.73205016915464682370 1.96122684512654421596  \\
  				5 1.73290875934971544403 2.96122647653791482014  \\
  				6 1.73376734954478517459 1.96122684512654421596  \\
  				7 2.59936313902457261094 2.46197022175260027410  \\
  				8 3.46495892850436026933 1.96122684512654488209  \\
  				9 3.46581751869942955580 2.96122647653791570832  \\
  				10 3.46667610889449795408 1.96122684512654554823  \\
  				11 4.33227189837428561248 2.46197022175259983001  \\
  				12 5.19786768785407282678 1.96122684512654288369  \\
  				13 5.19872627804914344551 2.96122647653791437605  \\
  				14 5.19958486824421228789 1.96122684512654421596  \\
  				15 6.06518065772399861402 2.46197022175259849774  \\
  				16 6.93077644720378494014 1.96122684512654221756  \\
  				17 6.93163503739885733523 2.96122647653791215561  \\
  				18 6.93249362759392528943 1.96122684512654155142  \\
  				19 7.79808941707371339191 2.46197022175259583321  \\
  				20 8.66368520655349882986 1.96122684512653866484  \\
  				21 8.66454379674857122495 2.96122647653790949107  \\
  				22 0.97032597650811958090 3.20302701275230772282  \\
  				23 0.27575758124369625923 3.92245369025308310285  \\
  				24 1.70869697454707125672 3.87742175565512958713  \\
  				25 1.92355445506011446533 2.90077634248431648345  \\
  				26 2.19931203630380833758 3.86200355619948565433  \\
  				27 2.89388043156823204782 3.14257687869871027431  \\
  				28 3.63225142960718327956 3.81697162160153213861  \\
  				29 3.84710891012022626612 2.84032620843071947903  \\
  				30 4.12286649136392036041 3.80155342214588864991  \\
  				31 4.81743488662834362657 3.08212674464511371397  \\
  				32 5.55580588466729530239 3.75652148754793557828  \\
  				33 5.77066336518033828895 2.77987607437712291869  \\
  				34 6.04642094642403193916 3.74110328809229208957  \\
  				35 6.74098934168845609349 3.02167661059151715364  \\
  				36 7.47936033972740688114 3.69607135349433857385  \\
  				37 7.69421782024044986770 2.71942594032352680244  \\
  				38 7.96997540148414262973 3.68065315403869641742  \\
  				39 7.79808941707371072738 1.46048346850048482715  \\
  				40 8.66454379674856767224 0.96122721371517128031  \\
  				41 6.93249362759392528943 1.96122684512653799871  \\
  				42 6.93163503739885733523 0.96122721371516584021  \\
  				43 6.93077644720378760468 1.96122684512653799871  \\
  				44 6.06518065772400039037 1.46048346850048393897  \\
  				45 5.19958486824421051153 1.96122684512653933098  \\
  				46 5.19872627804914078098 0.96122721371516450795  \\
  				47 5.19786768785407282678 1.96122684512653600031  \\
  				48 4.33227189837428650065 1.46048346850048216261  \\
  				49 3.46667610889449795408 1.96122684512654066324  \\
  				50 3.46581751869942955580 0.96122721371516717248  \\
  				51 3.46495892850435938115 1.96122684512653799871  \\
  				52 2.59936313902457261094 1.46048346850048460510  \\
  				53 1.73376734954478539663 1.96122684512653933098  \\
  				54 1.73290875934971366767 0.96122721371516961497  \\
  				55 1.73205016915464682370 1.96122684512654066324  \\
  				56 0.86645437967485894326 1.46048346850048615941  \\
  				57 0.00000000000000000000 0.96122721371517194644  \\
  				58 7.69421782024045253223 0.71942667750077782252  \\
  				59 8.38878621550487579839 0.00000000000000000000  \\
  				60 6.95584682220149908005 0.04503193459795417491  \\
  				61 6.74098934168845609349 1.02167734776876395486  \\
  				62 6.46523176044476333146 0.06045013405360009223  \\
  				63 5.77066336518034095349 0.77987681155436894276  \\
  				64 5.03229236714138750131 0.10548206865154810541  \\
  				65 4.81743488662834540293 1.08212748182236184746  \\
  				66 4.54167730538464997636 0.12090026810719266270  \\
  				67 3.84710891012022671021 0.84032694560796838967  \\
  				68 3.10873791208127459029 0.16593220270514649761  \\
  				69 2.89388043156823382418 1.14257761587595996211  \\
  				70 2.61812285032453928579 0.18135040216078968101  \\
  				71 1.92355445506011557555 0.90077707966156439490  \\
  				72 1.18518345702116500995 0.22638233675874316897  \\
  				73 0.97032597650812191237 1.20302774992955452404  \\
  				74 0.69456839526442792909 0.24180053621438463152  \\
  			};
  			
  			\addplot[no marks, 
  			nodes near coords={},
  			visualization depends on={value \thisrowno{0} \as \PunktI},
  			visualization depends on={value \thisrowno{1} \as \Scheitel},
  			visualization depends on={value \thisrowno{2} \as \PunktII},
  			visualization depends on={value \thisrowno{3} \as \Winkelradius},
  			visualization depends on={value \thisrowno{4} \as \Winkelfarbe},
  			visualization depends on={value \thisrowno{5} \as \Winkelname},
  			visualization depends on={value \thisrowno{6} \as \WinkelExzentrizitaet},
  			nodes near coords style={anchor=center,
  				path picture={
  					\begin{pgfonlayer}{bg}    
  					\draw pic [angle radius=\Winkelradius cm,%
  					fill=\Winkelfarbe!40, draw=\Winkelfarbe,
  					"$\Winkelname$", angle eccentricity =\WinkelExzentrizitaet,
  					text=\Winkelfarbe%
  					] {angle = P\PunktI--P\Scheitel--P\PunktII};
  					\end{pgfonlayer}
  			}},%
  			]
  			table[header=true, x expr=0.86645437967485894326, y expr=2.46197022175259894183]{
  				Punkt1 Scheitel Punkt2 Winkelradius[cm] Winkelfarbe Winkelname WinkelExz
  				3 1 4 0.5 blue!70!black {\gamma} 1.4 \\
  				6 7 8 0.5 blue!70!black \gamma{} 1.4 \\
  				10 11 12 0.5 blue!70!black {\gamma} 1.4 \\
  				14 15 16 0.5 blue!70!black {\gamma} 1.4 \\
  				18 19 20 0.5 blue!70!black {\gamma} 1.4 \\
  				1 2 22 0.4 green!50!black {\beta} 1.4 \\
  				24 22 23 0.4 orange!80!black {\alpha} 1.4 \\
  				28 27 26 0.4 orange!80!black {\alpha} 1.4 \\
  				32 31 30 0.4 orange!80!black {\alpha} 1.4 \\
  				36 35 34 0.4 orange!80!black {\alpha} 1.4 \\
  			};
  			
  			\end{axis}
  			
  			\draw[blue!70!white, line width=0.1, shorten >= -35, shorten <=-25] (P23) -- (P38);
  			\draw[blue!70!white, line width=0.1, shorten >= -25, shorten <=-40] (P74) -- (P59);
  			\draw[red!70!white, line width=0.1] (P3) -- (P20);
  			\draw[red!70!white, line width=0.1, shorten >= -65, shorten <=-65] (P11) -- (P48);
  			
  			\node[below, blue!70!white, font=\Large] at (P23) {$g_{1}\quad\quad$};
  			\node[above, blue!70!white, font=\Large] at (P74) {$g_{2}\quad\quad\quad\quad$};
  			\node[above, red!70!white, font=\Large] at (P30) {$g_{3}\;$};
  			
  			\node[above] at (P23) {$A$};
  			\node[left]  at (P3) {$B$};
  			\node[below] at (P74) {$C$};
  			\node[below] at (P59) {$D$};
  			\node[right] at (P20) {$E$};
  			\node[above] at (P38) {$F$};
  			\node[left]  at (P8) {$G\,\,$};
  			\node[right] at (P10) {$\,\,H$};
  			
  			\end{tikzpicture}
  		\captionof{figure}{The non-planar subgraph $G_4$.}
  	\end{minipage}
  \end{center}
  
  \newpage
  
  \begin{figure}[h]
  	\begin{subfigure}[b]{0.45\linewidth}
  	  \centering
  	  \begin{tikzpicture}
  		[y=0.5pt, x=0.5pt, yscale=-1.0, xscale=1.0]
  		\draw[line width=0.5pt]
  		(107.8208,945.8638) -- (114.4642,939.2204)
  		(70.1751,920.8753)arc(270.000:180.000:62.634)
  		(114.4642,939.2204)arc(315.000:270.000:62.634)
  		(70.1751,930.2704)arc(270.000:180.000:53.239)
  		(107.8208,945.8638)arc(315.000:270.000:53.239)
  		(114.4642,939.2204) -- (107.8208,945.8638)
  		(107.8208,945.8638)arc(135.000:90.000:62.634)
  		(114.4642,939.2204)arc(135.000:90.000:53.239)
  		(196.3990,945.8638) -- (189.7556,939.2204)
  		(296.6789,983.5095)arc(0.000:-90.000:62.634)
  		(234.0447,920.8753)arc(270.000:225.000:62.634)
  		(287.2838,983.5095)arc(0.000:-90.000:53.239)
  		(234.0447,930.2704)arc(270.000:225.000:53.239)
  		(189.7556,939.2204) -- (196.3990,945.8638)
  		(152.1099,964.2089)arc(90.000:45.000:62.634)
  		(152.1099,954.8138)arc(90.000:45.000:53.239)
  		(107.8208,1021.1552) -- (114.4642,1027.7986)
  		(7.5409,983.5095)arc(180.000:90.000:62.634)
  		(70.1751,1046.1438)arc(90.000:45.000:62.634)
  		(16.9360,983.5095)arc(180.000:90.000:53.239)
  		(70.1751,1036.7486)arc(90.000:45.000:53.239)
  		(114.4642,1027.7986) -- (107.8208,1021.1552)
  		(152.1099,1002.8101)arc(270.000:225.000:62.634)
  		(152.1099,1012.2052)arc(270.000:225.000:53.239)
  		(196.3990,1021.1552) -- (189.7556,1027.7986)
  		(234.0447,1046.1438)arc(90.000:-0.000:62.634)
  		(189.7556,1027.7986)arc(135.000:90.000:62.634)
  		(234.0447,1036.7486)arc(90.000:-0.000:53.239)
  		(196.3990,1021.1552)arc(135.000:90.000:53.239)
  		(189.7556,1027.7986) -- (196.3990,1021.1552)
  		(196.3990,1021.1552)arc(315.000:270.000:62.634)
  		(189.7556,1027.7986)arc(315.000:270.000:53.239);
  		\end{tikzpicture}
  	\end{subfigure}
  	\begin{subfigure}[b]{0.4\linewidth}
  	  \centering
  	    \begin{tikzpicture}
  	    [y=0.5pt, x=0.5pt, yscale=-1.0, xscale=1.0]
  	    \draw[line width=0.5pt]
  	    (71.5337,812.8510) -- (71.5337,803.4559)
  	    (71.5337,803.4559)arc(270.000:180.000:62.634)
  	    (71.5337,812.8510)arc(270.000:180.000:53.239)
  	    (124.7728,750.2168) -- (134.1679,750.2168)
  	    (71.5337,812.8510)arc(90.000:-0.000:62.634)
  	    (71.5337,803.4559)arc(90.000:-0.000:53.239)
  	    (187.4070,687.5826)arc(270.000:180.000:62.634)
  	    (187.4070,696.9777)arc(270.000:180.000:53.239)
  	    (303.2803,812.8510) -- (303.2803,803.4559)
  	    (365.9146,866.0901)arc(360.000:270.000:62.634)
  	    (356.5194,866.0901)arc(360.000:270.000:53.239)
  	    (250.0412,750.2168) -- (240.6461,750.2168)
  	    (240.6461,750.2168)arc(180.000:90.000:62.634)
  	    (250.0412,750.2168)arc(180.000:90.000:53.239)
  	    (250.0412,750.2168)arc(360.000:270.000:62.634)
  	    (240.6461,750.2168)arc(360.000:270.000:53.239)
  	    (71.5337,919.3292) -- (71.5337,928.7243)
  	    (8.8995,866.0901)arc(180.000:90.000:62.634)
  	    (18.2947,866.0901)arc(180.000:90.000:53.239)
  	    (124.7728,981.9634) -- (134.1679,981.9634)
  	    (134.1679,981.9634)arc(360.000:270.000:62.634)
  	    (124.7728,981.9634)arc(360.000:270.000:53.239)
  	    (124.7728,981.9634)arc(180.000:90.000:62.634)
  	    (134.1679,981.9634)arc(180.000:90.000:53.239)
  	    (303.2803,919.3292) -- (303.2803,928.7243)
  	    (303.2803,928.7243)arc(90.000:-0.000:62.634)
  	    (303.2803,919.3292)arc(90.000:-0.000:53.239)
  	    (250.0412,981.9634) -- (240.6461,981.9634)
  	    (303.2803,919.3292)arc(270.000:180.000:62.634)
  	    (303.2803,928.7243)arc(270.000:180.000:53.239)
  	    (187.4070,1044.5975)arc(90.000:0.000:62.634)
  	    (187.4070,1035.2024)arc(90.000:0.000:53.239);
  	    \end{tikzpicture}
  	\end{subfigure}
  	\caption{Non-planar examples of other forms besides ring graphs.}
  \end{figure}
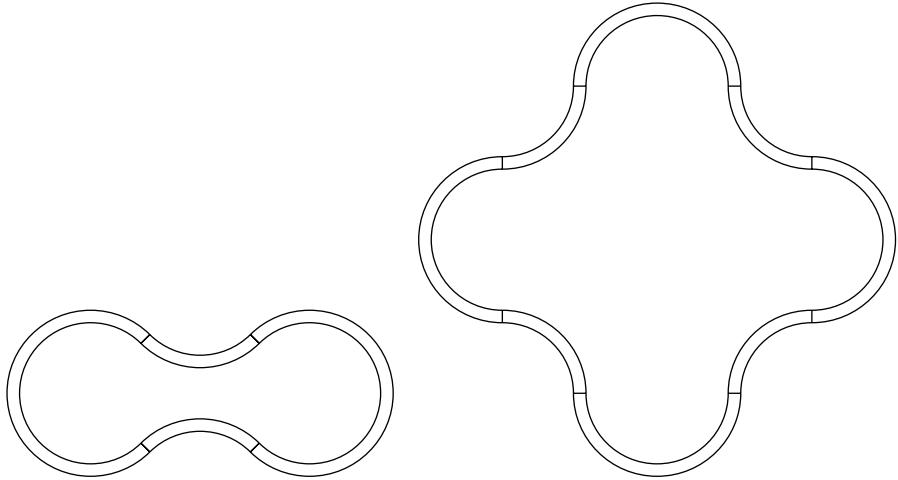
  
  Each graph from this article, as well as the ring graphs, can be viewed on the authors website \href{http://mikematics.de/matchstick_graphs_calculator.htm}{\textit{mikematics.de}}\footnote{http://mikematics.de/matchstick-graphs-calculator.htm}. The software\cite{MGC} we also used to verify the graphs runs directly in web browsers\footnote{For optimal functionality and design please use the Firefox web browser.}. The results of this article were first presented by the authors between September 17 and  November 3, 2018 in a graph theory internet forum\cite{Thread}.
  
  \section{\large{References}}
  
  \begingroup
  \renewcommand{\section}[2]{}
    
  \endgroup
  
\end{document}